\documentclass[twoside,leqno,12pt]{article}
\usepackage{ifthen,amsmath,amssymb,latexsym}
 \newtheorem{lemma}{Lemma}
 \newtheorem{proposition}{Proposition}

\numberwithin{equation}{section} \numberwithin{lemma}{section}
\numberwithin{theorem}{section} \numberwithin{corollary}{section}
\numberwithin{proposition}{section}
\numberwithin{definition}{section} \numberwithin{example}{section}
\numberwithin{remark}{section}
\def\TITLE       {About certain prime numbers}
\par
\smallskip
\smallskip
\def\LASTNAMEI   {\textbf{Savin}}
\def\FIRSTNAMEI  {\textbf{Diana}}

\smallskip

\def\ABSTRACT {
\begin{center}
\textbf{ABSTRACT}
\end{center}
We give a necessary condition for the existence of solutions of the Diophantine equation $p=x^{q}+ry^{q},$ with $p$, $q$, $r$ distinct odd prime natural numbers.
 
}
\def\CLASSIFICATION{MSC (2000): 11D41,11S15.}

\def\KEYWORDS{KEYWORDS: Diophantine equations; Kummer fields; cyclotomic fields.}

\begin{document}

\begin{center}
\par\TITLE
\par\FIRSTNAMEI \ \LASTNAMEI
\par
\par
Ovidius University, Constan\c ta, Romania

\end{center}
\par\
\par\ABSTRACT{}
\par\
\par\CLASSIFICATION{}
\par\KEYWORDS{}
\par\


 \section{Introduction}
 \label{Introduction}
  The representation of the prime natural numbers under the form
$x^{q}+ny^{q},$ for $q$$\in$ $\left\{2,3\right\}$ and $n$$\in$$\mathbb{N}^{*},$ has been already considered; we quote some papers in the Reference list.\\ 
   For example in D.Cox's book [2], there has been determined necessary and sufficient conditions for a prime $p$ to be written as
$p=x^{2}+ny^{2},$ where $n$$\in$$\mathbb{N}^{*}$.\\
 \begin{proposition} \label{pr1.1.}
 (\cite{2},pp.98,112) \textit{Let} $n$ \textit{be a square free positive integer which is not congruent to} $3$ \textit{modulo} $4.$
 \textit{Then there exists a monic irreducible polynomial} $f$$\in$$\mathbb{Z}[X]$ \textit{of degree} $h(\Delta),$ \textit{such that, if} $p$ \textit{is an odd prime not dividing} $n$ \textit{or the discriminant of} $f$ \textit{and} $E=HCF(K)$ \textit{is the Hilbert class field of} $K=\mathbb{Q}(\sqrt{-n}),$ \textit{the following statements are equivalent:}\\
(i) $p =x^{2}+ny^{2},$ \textit{for some} $x,y$$\in$$\mathbb{N}.$\\
(ii) $p$ \textit{completely splits in} $E.$\\
(iii) $\left(\frac{-n}{p}\right)=1$ \textit{and the congruence}
$ f(x)\equiv0$ (\textit{mod} $p$) \textit{has solutions in} $\mathbb{Z}.$\\
\textit{Moreover,} $f$ \textit{is the minimal polynomial of a real algebraic integer} $\alpha$ \textit{such that} $E=K(\alpha).$ 
\end{proposition}
\smallskip
In the above Proposition, $h(\Delta)$ denotes the number of classes of primitive positive definite forms of discriminant $\Delta,$ which is just the number of reduced forms (see D.Cox's book [2], p.29) and $\left(\frac{-n}{p}\right)$ denotes the Legendre symbol.\\
\smallskip
\ \ \ D.R. Heath-Brown \cite{1} considered the problem of finding primes of the form $x^{3}+2y^{3}$ and he obtained, by using techniques of analytic number theory, the following result:\\
\begin{proposition} \label{pr1.2.}
\textit{The number of primes of the form} $p=x^{3}+2y^{3}$ \textit{is infinite}.\\
\end{proposition}
D. Savin \cite{6} has proved the following result for $q=3:$
\begin{proposition} \label{pr1.3.}
\textit{Let} $K=\mathbb{Q}(\xi)$, $L=K(\sqrt[3]{r})$, \textit{where} $\xi$ \textit{is a primitive cubic root of unity and} $r$ \textit{is a prime number. If a prime number} $p$ \textit{can be written under the form} $p=x^{3}+ry^{3},$ \textit{then, for all} $P$$\in$Spec($\mathbb{Z}[\xi]$) \textit{above} $p$ (\textit{hence} $P$$\supseteq$$p\mathbb{Z}[\xi]$), \textit{the Artin symbol} $\left(\frac{L/K}{P}\right)$ \textit{is} $\textbf{1}_{L}.$\\
\end{proposition}
We shall define the Artin symbol of a prime ideal in Section 2.\\
\smallskip\\
In this paper, we generalize the above result for a prime $q\geq3$ (Lemma 3.1) and we use it for our main Proposition 3.1, in which we consider the general problem of determining whether or not there are primes $p$ of the form $p=x^{q}+ry^{q},$ with $q,r$ distinct odd prime numbers. We cannot solve the problem for any $q,r,$ but we find necessary conditions for the existence of solutions of the Diophantine equations $p=x^{q}+ry^{q},$ with $p,q,r$ distinct odd primes.These conditions (see Proposition 3.1) involve the ring of integers of a Kummer field $L=K(\sqrt[q]{r}),$ $K=\mathbb{Q}(\xi),$ where $\xi$ is a primitive root of unity of degree $q$ and $HCF(L)$ is its Hilbert class field.\\
First we recall some properties of integer rings of such fields in the next section.
\bigskip\\
\section{Notations and necessary concepts}
  The terminology we use here can be found in books of number theory, for example in \cite{3} or \cite{8}. \\
  Some properties of natural primes are studied by analysing the behavior of their ideals in the ring of integers in finite algebraic extensions of $\mathbb{Q.}$ We present only results which will be used in our paper; the proofs and other properties can be found in \cite{3} or \cite{8}.\\
  Denote by $K$ an algebraic extension of $\mathbb{Q}$ of degree $n$ and consider a prime number $p$ in $\mathbb{Z}.$ If $\mathbb{Z}_{K}$ denotes the ring of integers of $K$ over $\mathbb{Q},$ then $p\mathbb{Z}_{K}$ can be written as\\
  $$p\mathbb{Z}_{K}=\prod^{t}_{i=1}P^{e_{i}}_{i},$$ 
 where $P_{i},$ $i=1,2,...,t,$ are prime ideals in $\mathbb{Z}_{K}$ over $p,$ $e_{i},$ $i=1,2,...,t,$ are the ramification indices and $f_{i},$ $i=1,2,...,t,$ are the residual degrees (i.e. $f_{i}$ is the degree of the extension $\mathbb{Z}_{K}/P_{i}$ over $\mathbb{Z}/p\mathbb{Z}.$)\\
 If $K/\mathbb{Q}$ is a Galois extension, then the ramification indices $e_{i}$ are equal, $e_{i}=e,$ $i=1,2,...,t,$ as well as the residual degrees $f_{i},$ $f_{i}=f,$ $i=1,2,...,t,$ and $n=eft.$\\
\smallskip  
 If $K\subset L$ is an extension of algebraic number fields, $[L:K]=m,$ and there are no primes ramifying in $L,$ then $\Delta_{L}=\Delta^{m}_{K},$ where $\Delta_{L}$ and $\Delta_{K}$ are the absolute values of the discriminants of $L$ and $K,$ respectively.\\
 In D.Cox's book \cite{2}, the following proposition is proved.\\
 \begin{proposition} \label{pr2.1.}
(\cite{2}) \textit{Let} $K$ \textit{be an algebraic number field and} $P$$\in$ Spec($\mathbb{Z}_{K}$). \textit{Then} $P$ \textit{is totally split in the ring}
$\mathbb{Z}_{HCF(K)}$ \textit{if and only if} $P$ \textit{is a principal ideal in the ring} $\mathbb{Z}_{K}.$
\end{proposition}
For the cyclotomic field $\mathbb{Q}\left(\xi\right),$ where $\xi$ is a primitive root of unity of degree $r,$ $r\geq3,$ the following property of a prime number $p$ holds:\\
 \begin{proposition} \label{pr2.2.}(\cite{3})
\textit{Let} $r$$\in$$\mathbb{N},$ $r\geq3,$ $\xi$ \textit{be a primitive root of unity of order} $r$, $\mathbb{Z}\left[\xi\right]$ \textit{be the ring of integers of the cyclotomic field} $\mathbb{Q}\left(\xi\right)$. \textit{If} $p$ \textit{is a prime number which does not divide} $r$ \textit{and} $f$ \textit{is the smallest positive integer such that} $p^{f}$$\equiv$$1$$\pmod r$, \textit{then}\\
$$p\mathbb{Z}\left[\xi\right]=\prod^{t}_{i=1}P_{i},$$\\
 \textit{where} $t=\frac{\varphi\left(r\right)}{f}$ \textit{and} $P_{j},$ $j=1,2,...t,$
 \textit{are different prime ideals in the ring} $\mathbb{Z}\left[\xi\right]$.
\end{proposition}
 Now, if we consider an ideal $P$$\in$Spec($\mathbb{Z}\left[\xi\right]$) and $r$ is a prime natural number, then, in the ring of integers $A$ of the Kummer field $\mathbb{Q}\left(\sqrt[r]{\mu},\xi\right),$ with $\mu$$\in$$\mathbb{Z},$ the ideal $P$ is in one of the situations:\\ 
 (i)$P$ is the $r$-th power of a prime ideal, if the $r$-th power character $\left(\frac{\mu}{P}\right)_{r}=0;$\\
 (ii)$P$ decomposes in the product of $r$ different prime ideals of $A,$ if the $r$-th power character $\left(\frac{\mu}{P}\right)_{r}=1;$\\
 (iii)$P$ is a prime ideal of $A$, if the $r$-th power character $\left(\frac{\mu}{P}\right)_{r}$ is a root of order $r$ of unity, different from $1$.\\
\smallskip\\
Taking now in consideration a general finite Galois extension of algebraic number fields, $K\subset L,$ the following result is proved in D.Cox'book \cite{2}.  
\begin{proposition} \label{pr2.3.}
(\cite{2})\textit{Let} $K\subset L$ \textit{be a finite Galois extension of algebraic number fields and} $L=K(\alpha),$ $\alpha$$\in$$\mathbb{Z}_{L},$ $f$$\in$$\mathbb{Z}_{K}[X]$ \textit{be the minimal polynomial of} $\alpha$ \textit{over} $K$. \textit{If} $p$ \textit{is a prime element in} $\mathbb{Z}_{K}$
and $f$\textit{is separable modulo} $p$, \textit{then}:\\
(i) $p$ \textit{is unramified in the ring} $\mathbb{Z}_{L}.$\\
(ii)\textit{If} $f$$\equiv$$\prod^{t}_{i=1}f_{i}$$\pmod p$, \textit{where} $f_{i},$ $i=1,2,\ldots,t$, \textit{are pairwise distinct irreducible polynomials modulo} $p$, \textit{then} $P^{'}_{i}=p\mathbb{Z}_{L}+f_{i}(\alpha)\mathbb{Z}_{L}$ \textit{are prime ideals in the ring} $\mathbb{Z}_{L},$ $i=1,\ldots,t$, $P^{'}_{i}$$\neq$$P^{'}_{j},$ \textit{for} $i$$\neq$$j$ \textit{and} $p\mathbb{Z}_{L}=\prod^{t}_{i=1}P^{'}_{i}.$\\
\textit{Moreover, the polynomials }$f_{i},$ $i=1,\ldots,t$, \textit{have the same degree (equal to the residual degree} $f$).\\
(iii) $p$ \textit{is totally split in the ring} $\mathbb{Z}_{L}$ \textit{if and only if the congruence} $f(x)$ $\equiv$ $0$ $\pmod p$ \textit{has a solution in} $\mathbb{Z}_{K}.$
\end{proposition}
\begin{proposition} \label{pr2.4.}
(\cite{5})\textit{Let}$f$$\in$$\mathbb{Z}[X]$ \textit{be a monic irreducible polynomial and} $\alpha$ \textit{a root of} $f.$ \textit{We denote} $K=\mathbb{Q}\left(\alpha\right)$ \textit{and} ind($f$)$=[\mathbb{Z}_{K}:\mathbb{Z}[\alpha]]$ \textit{is the index of} $f.$ \textit{For a prime} $p$$\in$$\mathbb{N},$ \textit{we denote} $\overline{f}$$\in$$\mathbb{Z}_{p}[X]$ \textit{the polynomial obtained by reducing the coefficients of} $f$ \textit{modulo} $p.$ \textit{The prime decomposition of} $\overline{f}$ \textit{in} $\mathbb{Z}_{p}[X]$ \textit{is}
$\overline{f}=\varphi_{1}^{e_{1}}...\varphi_{r}^{e_{r}},$ \textit{where} $\varphi_{i}$$\in$$\mathbb{Z}_{p}[X],$ $i=1,...,t$ \textit{are distinct irreducible polynomials.}\\
\textit{If} $g_{i}$$\in$$\mathbb{Z}[X],$ $i=1,...,t$ \textit{are monic polynomials in} $\mathbb{Z}[X]$ \textit{such that} $\overline{g_{i}}=\varphi_{i},$ $i=1,...,t,$ \textit{and} $g=\frac{1}{p}\left(f-\prod^{t}_{i=1}g_{i}^{e_{i}}\right),$ \textit{then} $p$ \textit{does not divide} ind($f$) \textit{if and only if, for all} $i=1,2,...,t,$ \textit{either} $e_{i}=1$ \textit{or} $\varphi_{i}$ \textit{does not divide} $\overline{g}$ \textit{in} $\mathbb{Z}_{p}[X].$\\
\textit{It is known that, taking} $P_{i}=p\mathbb{Z}_{K}+g_{i}\left(\theta\right)\mathbb{Z}_{K},$ \textit{then} $P_{i}$$\in$Spec($\mathbb{Z}_{K}$),$i=1,...,t,$ and $p\mathbb{Z}_{K}=\prod^{t}_{i=1}P^{e_{i}}_{i}.$\\
\end{proposition}
It is well-known that, for an algebraic number field $K$ and a prime natural number $p,$ $p$ ramifies in $\mathbb{Z}_{K}$ if and only if $p$ divides the discriminant $\Delta_{K}.$\\
\smallskip\\
Let $L/K$ be a Galois extension and $P$$\in$Spec($\mathbb{Z}_{K}$) such that the extension $K$$\subset$$L$ is unramified in $P.$ Let $P^{'}$ be a prime ideal in the ring $\mathbb{Z}_{L}$ such that $P^{'}$ divides $P\mathbb{Z}_{L}.$ Then there exists a unique automorphism $\sigma$$\in$Gal$(L/K)$ such that\\
$$\sigma\left(x\right)\equiv x^{N(P)} (mod P^{'}),$$
where $N(P)=\left|\mathbb{Z}_{K}/P\right|$ is the norm of the ideal $P.$\\
\smallskip\\
The element $\sigma$ is called the \textbf{Artin symbol} of $P^{'}$ in the Galois extension $L/K$ and is denoted $\left(\frac{L/K}{P^{'}}\right).$\\
\smallskip\\
If the extension $K$$\subset$$L$ is Abelian, then $\left(\frac{L/K}{P^{'}}\right)$ does not depend on $P^{'}$ $\in$ Spec($\mathbb{Z}_{L}$) but only on $P=P^{'} \cap \mathbb{Z}_{K}$ and it is denoted $\left(\frac{L/K}{P}\right).$
The following result from \cite{7} will be useful in the proof of our result.
\begin{proposition} \label{pr2.5.}
\textit{If} $r,q$ \textit{are distinct prime  natural numbers},
 $\xi$ \textit{is a primitive root of order} $q$ \textit{of the unity,} $K=\mathbb{Q}\left(\xi\right)$ \textit{is the $q$-cyclotomic field and} $L=K(\sqrt[q]{r})$ \textit{is the Kummer field, then} 
 $$\left(\frac{L/K}{P}\right)\left(\sqrt[q]{r}\right)=\left(\frac{r}{P}\right)_{q}\sqrt[q]{r},$$
 \textit{for any} $P$$\in$Spec($\mathbb{Z}[\xi]$).
\end{proposition}
A transitivity property of ramification holds for Galois extensions $L$ and $M$ of the field $K,$ $K$$\subset$$M$$\subset$$L.$ A prime $p$ of $\mathbb{Z}_{K}$ splits completely in $L$ if and only if $p$ splits completely in $M$ and a prime ideal of $\mathbb{Z}_{M}$ containing $p$ splits completely in $L.$
\smallskip\\
Given the field finite extensions $\mathbb{Q}\subseteq$$k$$\subseteq$$K,$ there exists the following result on discriminants:
\begin{proposition} \label{pr2.6.}
(\cite{3}) \textit{Let} $K/\mathbb{Q}$ \textit{be a finite extension. Then, for every subfield} $k$$\subseteq$$K,$ \textit{we have} $\Delta_{k}^{[K:k]}$ \textit{divides} $\Delta_{K}.$\\
\end{proposition}
\bigskip
\section{Primes of the form $x^{q}+ry^{q}$} 
\smallskip
First we generalize our result from \cite{6} recalled in Proposition 1.3, for a prime $q\geq3.$
\begin{lemma} \label{lm3.1.}
\textit{Let} $p,q,r$ \textit{be distinct odd prime natural numbers}, $K=\mathbb{Q}(\xi)$, $L=K(\sqrt[q]{r})$, \textit{where} $\xi$ \textit{is a primitive root of unity of degree} $q$. \textit{If there exist integers} $x,y$ \textit{such that} $p=x^{q}+ry^{q},$ \textit{then the Artin symbol} $\left(\frac{L/K}{P}\right)=\textbf{1}_{L},$ \textit{for any} $P$$\in$Spec($\mathbb{Z}[\xi]$) \textit{which divides} $p\mathbb{Z}[\xi]$.
\end{lemma}
\textbf{Proof.} If there exist integers $x,y$ such that $p=x^{q}+ry^{q},$ then $x^{q}$$\equiv$$-ry^{q}$ (mod $p$). This implies that $x^{q}$$\equiv$$-ry^{q}$ (mod $P$), for any $P$$\in$Spec($\mathbb{Z}[\xi]$), $P|p\mathbb{Z}[\xi].$
We obtain for the $q$-power character that $\left(\frac{-ry^{q}}{P}\right)_{q}=1$. Since $\left(\frac{y^{q}}{P}\right)_{q}=1$ and $\left(\frac{-1}{P}\right)_{q}=1$, it results that $\left(\frac{r}{P}\right)_{q}=1.$\\
Using Proposition 2.4, we obtain that the Artin symbol, $\left(\frac{L/K}{P}\right),$ is equal to $\textbf{1}_{L},$ for all $P$$\in$Spec($\mathbb{Z}[\xi]$), $P|p\mathbb{Z}[\xi]$.\\
\smallskip\\
Now we state and prove our main result, establishing necessary conditions for the existence of primes of the form mentioned above.\\
\smallskip\\
In order to justify the hypotheses of the next Proposition, we give some examples of triples $(p,q,r)$ of primes satisfying the conditions: (i) $\overline{p}$ generates the group $(\mathbb{Z}^{*}_{q};\cdot);$ (ii) $(\exists)$ $x,y$$\in$$\mathbb{Z},$ such that $p=x^{q}+ry^{q}.$\\
\begin{enumerate}
\item [1)] $p=43,$ $q=5,$ $r=11.$ \\
Here $\overline{43}$$=\overline{3}$ generates the group $(\mathbb{Z}^{*}_{5},\cdot)$ and
$43=2^{5}+11\cdot1^{5}$. 
We obtain $\left(\frac{L/K}{43\mathbb{Z}[\xi]}\right)=\textbf{1}_{L},$ where
$~K=\mathbb{Q}(\xi),~L=\mathbb{Q}(\xi,\sqrt[5]{11}).$\\
\item [2)] $p=71,$ $q=3,$ $r=7.$ \\
Here $\overline{71}$$=\overline{2}$ generates the group $(\mathbb{Z}^{*}_{3},\cdot)$ and
$71=4^{3}+7\cdot1^{3}$.\\
We obtain $\left(\frac{L/K}{71\mathbb{Z}[\epsilon]}\right)=\textbf{1}_{L},$ where
$~K=\mathbb{Q}(\epsilon),~L=\mathbb{Q}(\epsilon,\sqrt[3]{7}).$\\
\item [3)] $p=131,$ $q=3,$ $r=13.$ \\
Here $\overline{131}$$=\overline{2}$ generates the group $(\mathbb{Z}^{*}_{3},\cdot)$ and
$131=3^{3}+13\cdot2^{3}$. 
We obtain $\left(\frac{L/K}{131\mathbb{Z}[\epsilon]}\right)=\textbf{1}_{L},$ where
$~K=\mathbb{Q}(\epsilon),~L=\textbf{Q}(\epsilon,\sqrt[3]{13}).$\\
\item [4)] $p=197,$ $q=3,$ $r=7.$ \\
Here $\overline{197}$$=\overline{2}$ generates the group $(\mathbb{Z}^{*}_{3},\cdot)$ and
$197=2^{3}+7\cdot3^{3}$. \\
We obtain $\left(\frac{L/K}{71\textit{Z}[\epsilon]}\right)=\textbf{1}_{L},$ where
$~K=\mathbb{Q}(\epsilon),~L=\mathbb{Q}(\epsilon,\sqrt[3]{7}).$\\
\item [5)] $p=353,$ $q=5,$ $r=11.$ \\
Here $\overline{353}$$=\overline{3}$ generates the group $(\mathbb{Z}^{*}_{5},\cdot)$ and
$353=1^{5}+11\cdot2^{5}$. 
We obtain $\left(\frac{L/K}{353\mathbb{Z}[\xi]}\right)=\textbf{1}_{L},$ where
$~K=\mathbb{Q}(\xi),~L=\mathbb{Q}(\xi,\sqrt[5]{11}).$\\
\item [6)] $p=16903,$ $q=5,$ $r=3.$ \\
Here $\overline{16903}$$=\overline{3}$ generates the group $(\mathbb{Z}^{*}_{5},\cdot)$ and
$16903=7^{5}+3\cdot2^{5}$. 
We obtain $\left(\frac{L/K}{16903\textit{Z}[\xi]}\right)=\textbf{1}_{L},$ where
$~K=\mathbb{Q}(\xi),~L=\mathbb{Q}(\xi,\sqrt[5]{3}).$\\
\item [7)] $p=127277,$ $q=7,$ $r=3.$ \\
Here $\overline{127277}$$=\overline{3}$ generates the group $(\mathbb{Z}^{*}_{7},\cdot)$ and
$127277=5^{7}+3\cdot4^{7}$. 
We obtain $\left(\frac{L/K}{127277\mathbb{Z}[\xi]}\right)=\textbf{1}_{L},$ where
$~K=\mathbb{Q}(\xi),~L=\mathbb{Q}(\xi,\sqrt[7]{3}).$\\
\item [8)] $p=187387,$ $q=11,$ $r=5.$ \\
Here $\overline{187387}$$=\overline{2}$ generates the group $(\mathbb{Z}^{*}_{11},\cdot)$ and
$187387=3^{11}+5\cdot2^{11}$. 
We obtain $\left(\frac{L/K}{187387\mathbb{Z}[\xi]}\right)=\textbf{1}_{L},$ where
$~K=\mathbb{Q}(\xi),~L=\mathbb{Q}(\xi,\sqrt[11]{5}).$\\
\end{enumerate}
A simply computation in MATHEMATICA, using for example the procedure\\ Do[Print[FactorInteger[$i^{11}+5*j^{11}$]], $\left\{i,10\right\},$ $\left\{j,10\right\}$], gives more and more examples of primes satisfying the conditions (i) and (ii).\\
We give the precise results (obtained using the program MAGMA) for the $4$-th example:\\
$\alpha=65$ is a solution for the congruence $\alpha^{3}$$\equiv$$7$ (mod $197$), so the decomposition of the ideal $197A$ in a product of prime ideals of $A$ is:\\
$$pA=(197,65-\sqrt[3]{7})\cdot(197,65-\xi\sqrt[3]{7})\cdot(197,65-\xi^{2}\sqrt[3]{7}),$$
where $A=\mathbb{Z}_{L}.$\\
After this we obtained with MAGMA that the ideals $P_{i}=(197,65-\xi^{i}\sqrt[3]{7}),$ $i=0,1,2$ are totally split in the ring of integers of the Hilbert field $HCF(L),$ so $p=65$ which is prime in $\mathbb{Z}_{K},$ is totally split in $\mathbb{Z}_{HCF(L)}.$\\
For all these examples, with a simple verification with MAGMA, we get that the ideal $p\mathbb{Z}_{HCF(L)}$ belongs to $S_{HCF(L)/K}.$ 
\begin{proposition} \label{pr3.1.}
\textit{Let} $p,q,r$ \textit{be distinct odd prime natural numbers, such that} $\overline{p}$ \textit{generates the group} $(\mathbb{Z}^{*}_{q};\cdot),$ $K=\mathbb{Q}(\xi)$, $L=K(\sqrt[q]{r})$, \textit{where} $\xi$ \textit{is a primitive root of unity of degree} $q$. \textit{Let} $P_{K}$ \textit{be the set of all finite primes of} $K$ \textit{and} $S_{HCF(L)/K}=$ $P \in P_{K}: P$ \textit{splits completely in} $HCF(L)\}.$ \textit{If there exist rational integers} $x,y$ \textit{such that} $p=x^{q}+ry^{q},$ \textit{then:}\\ i)
\textit{The ideal generated by} $p$ belongs to $S_{HCF(L)/K}.$\\
ii) \textit{There exists a real algebraic integer} $\alpha$ \textit{such that} $HCF(L)=L(\alpha)$ \textit{and if we denote with} $f$ \textit{the minimal polynomial of } $\alpha$ \textit{over} $L,$ \textit{then the congruence 
} $f(x)\equiv0$ (\textit{mod} $p$) \textit{has solutions in} $\mathbb{Z}[\xi].$
\end{proposition} 
\textbf{Proof.} i) We denote by $A$ the ring of integers of the field $L$.\\
Since $<\overline{p}>$$=$$\mathbb{Z}^{*}_{q},$ it results that $p$ remains prime in the ring $\mathbb{Z}[\xi].$\\
We have\\
$$pA=(x+y\sqrt[q]{r})A(x+\xi y\sqrt[q]{r})A\cdots(x+\xi^{q-1}y\sqrt[q]{r})A.$$\\
But, according to Lemma 3.1, $\left(\frac{r}{pZ[\xi]}\right)_{q}=1$ and we obtain that $P_{i}=(x+\xi^{i}y\sqrt[q]{r})A$$\in$Spec($A$), for all $i=1,...,q-1.$ Using Proposition 2.1 we get that $P_{i}=(x+\xi^{i}y\sqrt[q]{r})A$$\in$Spec($A$) is totally split in the ring $\mathbb{Z}_{HCF(L)},$ so (according to a result from Section 2) $p$ is totally split in the ring $\mathbb{Z}_{HCF(L)}.$\\ 
ii) First we show that $f$ is separable modulo $p.$\\
From i), we obtain that $p$ does not ramify in $A.$
This is equivalent with the fact that $p$ does not divide $\Delta_{L}.$\\
Since $L$ $\subset$ $HCF(L)$ is an Abelian unramified extension, we obtain that $\Delta_{HCF(L)}=$$\Delta_{L}^{h_{L}}.$\\
Since $p$ does not divide $\Delta_{L},$ then $p$ does not divide $\Delta_{HCF(L)}.$ \\
It is known that $\mathbb{Q}$$\subset$$HCF(L)$ is a Galois extension. Let $G$ be the Galois group of the extension $HCF(L)/\mathbb{Q}$ and let $\tau$ denote the complex conjugation. The $H=\left\{\textbf{1},\tau\right\}$ is a subgroup in $G.$\\
Let $F=HCF(L)^{H}=\left\{x \in HCF(L):u(x)=x,\forall u\in H \right\}.$ It is results \\
$F=\left\{x \in HCF(L):\tau(x)=x\right\},$ so $F$ is real and $[HCF(L):F]=\left|H\right|=2.$\\
Since $L$ is not real and $F$ is real, we obtain that there exists a real algebraic integer $\alpha$ such that
$HCF(L)=L(\alpha),$ $F=\mathbb{Q}(\alpha).$\\
We denote with ind($f$)$=[\mathbb{Z}_{F}:\mathbb{Z}[\alpha]].$ Using i) and Proposition 2.4, it results $p$ does not divide ind($f$).\\
Using Proposition 2.6 we obtain that $\Delta_{F}$ divides $\Delta_{HCF(L)},$ so $p$ does not divide $\Delta_{F}.$\\
Knowing that disc($f$)=(ind $(f))^{2}\cdot\Delta_{F},$ we obtain that $p$ does not divide disc($f$), so $f$ is separable modulo $p.$ Applying Proposition 2.3, we obtain that the congruence $f(x)\equiv0$ (mod $p$) has a solution in A.\\
Since $\left(\frac{r}{pZ[\xi]}\right)_{q}=1,$ it results that the number of prime ideals factorizing $pA$ is $q,$ so the residual degrees $f_{P_{i}/pA}$ are equal to $1.$ This implies the field isomorphism $A/P_{i}$$\cong$$\mathbb{Z}[\xi]/p\mathbb{Z}[\xi].$ Therefore the fact that the congruence $f(x)\equiv0$ (mod $p$) has a solution in $A$ is equivalent to the fact that the congruence $f(x)\equiv0$ (mod $p$) has a solution in $\mathbb{Z}[\xi]$. \\
\smallskip\\
  
\bigskip
{Faculty of Mathematics and Computer Science ,\\
  Department of Mathematics\\
"Ovidius" University of Constanta}\\
Bd. Mamaia 124, Constanta, 900527\\
Romania\\
{e-mails: Savin.Diana@univ-ovidius.ro\\
\ \ \ \ \ } 
  \bigskip\\

\begin{thebibliography}{l} 
\bibitem{1} D.R. Heath-Brown, \textit{Primes Represented by} $x^{3}+ 2y^{3}$, Acta Math.,\textbf{186}(2001),1-84.
\bibitem{2} D. Cox, \textit{Primes of the Form} $x^{2}+ny^{2}$, A Wiley - Interscience Publication, New York, 1989.
\bibitem{3} F. Lemmermeyer, \textit{Reciprocity laws, from Euler to Eisenstein }, Springer-Verlag, Heidelberg, 2000.
 \bibitem{4}M. Magioladitis, \textit{Primes of the form} $x^{2}+ny^{2}$, Arbeitsgemeinschaft zur Klassenkorpertheorie, University of Duisburg-Essen, 2004.
\bibitem{5} J. Montes, E. Nart, \textit{On a theorem of Ore}, J. Algebra, \textbf{146} (1992), 318-334.
\bibitem{6} D. Savin, \textit{Artin symbol of the Kummer fields}, J. Creative Math. Informatics,\textbf{16} (2007), 63-69.
\bibitem{7} P. Stevenhagen, \textit{Kummer Theory and Reciprocity Laws}, Universiteit Leiden, 2005.
\bibitem{8} L.C. Washington, \textit{Introduction to Cyclotomic Fields}, Springer, New York, 1997.
\end{thebibliography}
 \end{document}